\documentclass[12pt]{ amsart}
\usepackage{amsmath}
\usepackage{amsfonts}
\usepackage{amssymb}
\usepackage{amsthm}
\usepackage{enumitem}

\newtheorem{theorem}{Theorem}

\newtheorem*{rem}{Remark}

\begin{document}
\title[On Diophantine triples containing a triangular number]{On Diophantine triples containing a triangular number}
\author{Marija {Bliznac Trebje\v{s}anin}}

\begin{abstract}
A general construction yielding infinitely many families of $D(m^2)$-triples of triangular numbers is presented. Moreover, each triple obtained from this construction contains the same triangular number $T_n$.
\end{abstract}

\maketitle

\noindent 2020 {\it Mathematics Subject Classification:} 11B37, 11D09
\\ \noindent Keywords: triangular numbers, recurrences, Diophantine m-tuples

\section{Introduction}

Let $n\neq0$ be an integer. We call the set of $m$ distinct positive integers a $D(n)$-$m$-tuple (or a Diophantine $m$-tuple with the property $D(n)$) if the product of any two of its distinct elements increased by $n$ is a perfect square. 
When $n=1$, the so-called classical case introduced by Diophantus, such a set is simply referred to as a Diophantine m-tuple.

The first example of a Diophantine quadruple was found by Fermat, and it was the set $\{1, 3, 8, 120\}$. One of the most studied questions is how large those sets can be. He, Togb\'{e} and Ziegler in \cite{htz} proved the nonexistence of $D(1)$-quintuples, while the author and Filipin proved the nonexistence of $D(4)$-quintuples in \cite{petorke_btf}. In the case $n=-1$, the nonexistence of quadruples was established in \cite{bcm}. 
Another well-studied variant of the problem is to determine whether elements of certain interesting sequences can occur in $D(n)$-$m$-tuples. Dujella in \cite{duje_gen} listed some quadruples with the properties $D(1)$, $D(4)$, $D(9)$ and $D(64)$ that contain two or three Fibonacci numbers. In \cite{fujita_luca_fib} Fujita and Luca showed that there is no Diophantine quadruple consisting solely of Fibonacci numbers. Numerous other variants and generalizations of the concept of Diophantine m-tuples have been investigated; for an overview of these problems and related references, see Section 1.5 of \cite{duje_book}.

Let $n$ be a positive integer. The triangular number $T_n$ is a figurate number that can be represented as a triangular grid of points with 
$n$ rows, where the first row contains a single element and each following row contains one more element than the previous one. Thus, a triangular number is a sum of all positive integers less than or equal to a given positive integer $n$, which implies the explicit expression
$$T_n=\frac{n(n+1)}{2}.$$

From the main result of \cite{mnde}, it follows that  
$$\{8,T_n,T_{n+4},T_{4n^2+20n+8}\}$$  
is a Diophantine quadruple for every positive integer $n$, which implies that  
$$\{T_n,T_{n+4},T_{4n^2+20n+8}\}$$  
is a Diophantine triple of triangular numbers.
These Diophantine triples are also studied in \cite{hamtat}.

It is easy to verify that for each positive integer $m$, triple 
$$\{T_n,T_{n+m},T_{4(n^2+(4m+1)n+2m)}\}$$
is a $D(m^2)$-triple, for any positive integer $n$. 

In the rest of the paper, we present a construction yielding infinitely many 
$D(m^2)$-triples of triangular numbers. Moreover, each triple contains the initial triangular number 
$T_n$ thereby establishing the following main result.

\begin{theorem}\label{tm_1}
Let $m$ be a positive integer. For every positive integer $n$, triangular number $T_n$ is a member of infinitely many $D(m^2)$-triples of triangular numbers.
\end{theorem}


\section{Constructing Diophantine triples}\label{sec_prvi}

Let $m$ be a positive integer. For every positive integer $n$, pair $\{T_n,T_{n+4m}\}$ is a $D(m^2)$-pair since 
$$T_n\cdot T_{n+4m}+m^2=\left(\frac{n^2+(4m+1)n+2m}{2}\right)^2=:r_1^2.$$
We can extend this $D(m^2)$-pair to a $D(m^2)$-triple $\{T_n,T_{n+4m},T_{8r_1}\}$. 
Then 
\begin{align*}
     T_nT_{8r_1}+m^2&=(2r_1(2n+1)-m)^2=s_1^2,  \\
    T_{n+4m}T_{8r_1}+m^2&=(2r_1(2(n+4m)+1)+m)^2=t_1^2.
\end{align*}

Let $N_0=n$, $N_1=n+4m$, $N_2=8r_1$. Define $N_3=8s_1+N_1$. Then the triple $\{T_{N_0},T_{N_2},T_{N_3}\}$ is a $D(m^2)$-triple. It holds
\begin{align*}
T_{N_0}T_{N_2}+m^2&=s_1^2=:r_2^2,\\
     T_{N_0}T_{N_3}+m^2&=(2s_1(2N_0+1)-r_1)^2=:s_2^2,  \\
    T_{N_2}T_{N_3}+m^2&=(2s_1(2N_2+1)+t_1)^2=:t_2^2.
\end{align*}
Let's define recursively
$$
    N_{k+2}=8s_k+N_k,\quad k\geq 1.
$$
We will prove that  
$$\{T_{N_0},T_{N_{k+1}},T_{N_{k+2}}\}$$
is $D(m^2)$-triple with 
\begin{align}
    s_{k+1}&=2s_k(2N_0+1)-s_{k-1},\label{rek_s}\\
    t_{k+1}&=2s_k(2N_{k+1}+1)+t_k.\label{rek_t}
\end{align}

Recurrence \eqref{rek_s} can be solved to obtain an explicit formula for $s_k$. Note that we can take $s_0:=r_1$ and $s_{-1}:=m$. Then, it holds 
$$s_k=\lambda_1\alpha^{k+1}+\lambda_2\beta^{k+1},\ k\geq 1,$$
where 
\begin{align*}
    &\alpha=2N_0+1+2\sqrt{N_0(N_0+1)},\quad \lambda_1=\frac{m}{2}+\frac{1}{8}\sqrt{N_0(N_0+1)},\\
    &\beta=2N_0+1-2\sqrt{N_0(N_0+1)},\quad \lambda_2=\frac{m}{2}-\frac{1}{8}\sqrt{N_0(N_0+1)}.\\
\end{align*}
It also holds $\alpha=(\sqrt{N_0+1}+\sqrt{N_0})^2$ and $\beta=(\sqrt{N_0+1}-\sqrt{N_0})^2$, hence  $\beta=\alpha^{-1}$ and $\beta/(\beta^2-1)=-\alpha/(\alpha^2-1)$.

By carefully considering the cases when $k$ is even or odd, we obtain
$$N_k=\frac{8\alpha}{\alpha^2-1}(\lambda_1\alpha^k-\lambda_2\beta^k)-\frac{1}{2},$$
and, similarly,
$$t_k=\frac{32\alpha^3\lambda_1^2}{(\alpha^2-1)^2}\alpha^{2k}+\frac{32\alpha\lambda_2^2}{(\alpha^2-1)^2}\beta^{2k}-\frac{32\alpha}{(\alpha^2-1)^2}(\alpha^2\lambda_1^2+\lambda_2^2)+m.$$

Suppose that $\{T_{N_0},T_{N_k},T_{N_{k+1}}\}$ (and all its preceding triples) forms a $D(m^2)$-triple, and that \eqref{rek_s} and \eqref{rek_t} hold for this triple.  
 To show that the triple $\{T_{N_0},T_{N_{k+1}},T_{N_{k+2}}\}$ also has the property $D(m^2)$ we will use the previously obtained expressions.  
We need to show that
\begin{equation}\label{eq_s}
T_{N_0}T_{N_{k+2}}+m^2=s_{k+1}^2
\end{equation}
and
\begin{equation}\label{eq_t}
T_{N_{k+1}}T_{N_{k+2}}+m^2=t_{k+1}^2
\end{equation}
hold,
where $s_{k+1}$ and $t_{k+1}$ are given by \eqref{rek_s} and \eqref{rek_t}. 

The following simple properties of triangular numbers will be useful in simplifying the proof:
$$T_{8x+y}=32x^2+4x(2y+1)+T_{y}$$
and
$$8T_n+1=(2N+1)^2.$$
Let us first consider \eqref{eq_s}. 
By the definition of the $D(m^2)$ property and the fact that $T_{N_0}T_{N_k}+m^2=s_{k-1}^2$, we have 
\begin{align*}
T_{N_0}T_{N_{k+2}}+m^2&=T_{N_0}T_{8s_k+N_k}+m^2\\
&=32s_k^2T_{N_0}+4s_k(2N_k+1)T_{N_0}+s_{k-1}^2.
\end{align*}
On the other hand, from \eqref{rek_s} we have 
\begin{align*}
s_{k+1}^2&=4s_k^2(2N_0+1)^2-4s_ks_{k-1}(2N_0+1)+s_{k-1}^2,\\
&=32s_k^2T_{N_0}+4s_k^2-4s_ks_{k-1}(2N_0+1)+s_{k-1}^2.
\end{align*}
Hence, the proof reduces to showing that
\begin{equation}\label{jedn_s}
    (2N_k+1)T_{N_0}=s_k-s_{k-1}(2N_0+1) 
\end{equation}
holds. This can be easily verified by substituting the explicit expressions for $N_k$, $s_k$ and $s_{k-1}$.

Similarly, proving \eqref{eq_t} reduces to showing that
\begin{equation}
T_{N_{k+1}} (2N_k + 1) = s_k + t_k (2N_{k+1} + 1)
\end{equation}
holds. Once again, by substituting the explicit expressions for the involved quantities, we verify that this equality is satisfied.

This completes the proof of Theorem \ref{tm_1}.

\begin{rem}
   Obviously, for some positive integers $m$ and $n$, these are not the only $D(m^2)$-triples in triangular numbers containing $T_n$.  For instance, $T_{n+4m}$ appears as the second element in the $D(m^2)$-triple $\{T_n,T_{n+4m},T_{8r}\}$, which is not part of the sequence of $D(m^2)$-triples produced by the construction starting with $T_{N_0} = T_{n+4m}$. 
But there are examples where $T_n$ is the smallest element in a triple that does not arise from the construction. For example, the Diophantine triples $\{T_1,T_{15},T_{90}\}$ and $\{T_2,T_{15},T_{153}\}$, as well as the $D(9)$-triple $\{T_1,T_{63},$ $T_{370}\}$, all exhibit this property.

\end{rem}

\bigskip
\textbf{Acknowledgement:} This work was supported by the Croatian Science Foundation grant no. IP-2022-10-5008.\\

\bigskip
Marija Bliznac Trebje\v{s}anin\\
Faculty of Science, University of Split,\\ Ru\dj{}era Bo\v{s}kovi\'{c}a 33, 21000 Split, Croatia \\
Email: marbli@pmfst.hr \\
\\


\begin{thebibliography}{99}

\bibitem{petorke_btf} M.~Bliznac Trebje\v{s}anin, A.~Filipin, {\it Nonexistence of $D(4)$-quintuples}, {  J.~Number Theory}, \textbf{194} (2019), 170--217.

\bibitem{bcm} N.~C.~Bonciocat, M.~Cipu and M.~Mignotte, {\it There is no Diophantine D(-1)-quadruple}, J. London Math. Soc. \textbf{105} (2022), 63--99.

\bibitem{mnde}M.~N.~Deshpande, {\it One property of triangular numbers}, Portugaliae Math. \textbf{55} (1998), 381--383.

\bibitem{duje_gen} A. Dujella, {\it Generalization of a problem of Diophantus}, Acta Arith. \textbf{65} (1993), 15--27.

\bibitem{duje_book} A. Dujella, {\it Diophantine $m$-tuples and Elliptic Curves}, Springer, Cham, 2024.

\bibitem{fujita_luca_fib} Y.~Fujita, F.~Luca, {\it There are no Diophantine quadruples of Fibonacci numbers}, Acta Arith. \textbf{185} (2018), 19--38.

\bibitem{hamtat} A.~Hamtat, {\it Diophantine triples in triangular numbers}, preprint, 2025.

\bibitem{htz} B.~He, A.~Togb\'{e}, V.~Ziegler, {\it There is no Diophantine quintuple}, {  Trans.~Amer.~Math.~Soc.~}\textbf{371} (2019), 6665--6709.


\end{thebibliography}
\end{document}